\documentclass[12pt]{article}
\usepackage{amsfonts}
\usepackage{latexsym,amssymb}
\usepackage{mathrsfs}
\usepackage{graphicx}
\usepackage{psfrag}
\usepackage{amsmath, amsbsy}
\usepackage{amsopn, amstext}
\usepackage{multicol}
\usepackage{latexsym,amssymb}

\allowdisplaybreaks

\usepackage{graphicx}
\usepackage{hyperref}

\def\sqr#1#2{{\vcenter{\vbox{\hrule height.#2pt
              \hbox{\vrule width.#2pt height#1pt \kern#1pt \vrule width.#2pt}
              \hrule height.#2pt}}}}
\def\signed #1{{\unskip\nobreak\hfil\penalty50
              \hskip2em\hbox{}\nobreak\hfil#1
              \parfillskip=0pt \finalhyphendemerits=0 \par}}
\def\endpf{\signed {$\sqr69$}}

\def\dbR{{\mathop{\rm l\negthinspace R}}}

\def\3n{\negthinspace \negthinspace \negthinspace }
\def\2n{\negthinspace \negthinspace }
\def\1n{\negthinspace }

\def\ds{\displaystyle}
\def\dbA{{\mathbb{A}}}

\def\dbN{{\mathbb{N}}}

\def\dbR{{\mathbb{R}}}

\def\={\buildrel \triangle \over =}

\def\mE{\mathbb{E}}
\def\mR{\mathbb{R}}
\def\mN{\mathbb{N}}
\def\l{\lambda}

 \def\n{\nabla}
\def\si{\sigma}
\def\t{\times}

\def\ns{\noalign{\ss} }

\def\G{\Gamma}

\def\O{\Omega}

\def\cA{{\cal A}}

\def\cC{{\cal C}}

\def\cF{{\cal F}}

\def\cL{{\cal L}}

\def\cS{{\cal S}}

\def\cV{{\cal V}}

\def\cX{{\cal X}}
\def\cY{{\cal Y}}

\def\mE{{\mathbb{E}}}

\def\no{\noindent}

\def\ms{\medskip}
\def\bs{\bigskip}
\def\q{\quad}
\def\qq{\qquad}
\def\hb{\hbox}

\def\max{\mathop{\rm max}}

\def\sup{\mathop{\rm sup}}

\def\cd{\cdot}
\def\cds{\cdots}

\def\as{\hbox{\rm a.s.{ }}}

\def\|{\Big |}
\def\({\Big (}
\def\){\Big )}
\def\[{\Big[}
\def\]{\Big]}
\def\be{\begin{equation}}
\def\bel{\begin{equation}\label}
\def\ee{\end{equation}}
\def\bt{\begin{theorem}}
\def\bcd{\begin{condition}}
\def\ecd{\end{condition}}
\def\et{\end{theorem}}
\def\bc{\begin{corollary}}
\def\ec{\end{corollary}}
\def\bde{\begin{definition}}
\def\ede{\end{definition}}
\def\bl{\begin{lemma}}
\def\el{\end{lemma}}
\def\bp{\begin{proposition}}
\def\ep{\end{proposition}}
\def\br{\begin{remark}}
\def\er{\end{remark}}
\def\ba{\begin{array}}
\def\ea{\end{array}}
\def\ed{\end{document}}
\def\ns{\noalign{\ms}}
\def\ds{\displaystyle}

\def\square#1{\vbox{\hrule\hbox{\vrule height#1%
     \kern#1\vrule}\hrule}}
\def\rectangle#1#2{\vbox{\hrule\hbox{\vrule height#1%
     \kern#2\vrule}\hrule}}

\font\tenbb=msbm10 \font\sevenbb=msbm7 \font\fivebb=msbm5

\newfam\bbfam
\scriptscriptfont\bbfam=\fivebb \textfont\bbfam=\tenbb
\scriptfont\bbfam=\sevenbb

\newtheorem{lemma}{Lemma}[section]
\newtheorem{remark}{Remark}[section]

\newtheorem{theorem}{Theorem}[section]
\newtheorem{corollary}{Corollary}[section]

\newtheorem{definition}{Definition}[section]
\newtheorem{proposition}{Proposition}[section]
\newtheorem{condition}{Condition}[section]
\makeatletter
   
   \@addtoreset{equation}{section}
\makeatother

\begin{document}
\title{\bf Exact Controllability for  Stochastic Transport Equations\thanks{This work
is partially supported by the NSF of China
under grant 11101070 and the Fundamental
Research Funds for the Central Universities in
China under grants ZYGX2012J115.  \ms} \ms}

\author{Qi L\"u\thanks{School of Mathematical Sciences, University of Electronic
 Science and Technology of China, Chengdu, 610054, China.
  {\small\it E-mail:} {\small\tt luqi59@163.com}.} }

\date{}

\maketitle

\begin{abstract}\no
This paper is addressed to studying the exact
controllability for stochastic transport
equations  by two controls: one is a boundary
control imposed on the drift term and the other
is an internal control imposed on the diffusion
term.
 By means of the
duality argument, this controllability problem
can be reduced to an observability problem for
backward stochastic transport equations, and the
desired observability estimate is obtained by a
new global Carleman estimate. Also, we present
some results about the lack of exact
controllability, which show that the action of
two controls is necessary. To some extent, this
indicates that the controllability problems for
stochastic PDEs differ from their deterministic
counterpart.
\end{abstract}

\bs

\no{\bf 2010 Mathematics Subject
Classification}.  Primary 93B05; Secondary
93B07, 93E20, 60H15.

\bs

\no{\bf Key Words}. Stochastic transport
equations, exact controllability,
observability, Carleman estimate.

\ms

\section{Introduction}

Let $T>0$ and $G\subset \mR^d$\;($d\in\mN$) be a
strictly convex bounded domain with a $C^1$
boundary $\G$. Denote by
$\nu(x)=(\nu^1(x),\cds,\nu^d(x))$ the unit
outward normal vector of $G$ at $x\in\G$. Let
$\bar x_1,\bar x_2\in \G$  satisfy  that
$$|\bar x_1 - \bar
x_2|_{\mR^d}=\ds\max_{ x_1,x_2\in \overline
G}|x_1-x_2|_{\mR^d}.$$ Without loss of
generality, we assume that $0\in G$ and $0=\bar
x_1 + \bar x_2$.   Put $R=\ds\max_{x\in \G}
|x|_{\mR^d}$. Let
$$
S^{d-1}\=\{x\in\mR^d:\,|x|_{\mR^d}=1\}.
$$
Denote by
$$
\G_{S}^- = \{ (x,U)\in\G\times S^{d-1}:\,
U\cd\nu(x)\leq 0 \},\;\; \G_{S}^+=\big(\G\times
S^{d-1}\big)\setminus \G_{S}^-.
$$
Let us define a Hilbert space $
L^2_{w}(\G_{S}^-) $ as the completion of all
$h\in C_0^\infty(\G_S^-\times S^{d-1})$ with the
norm
$$
|h|_{L^2_{w}(\G_{S}^-)}\=\(\ds-
\int_{\G_{S}^-}U\cd\nu |h|^2
d\G_{S}^-\)^{\frac{1}{2}},
$$
where $d\G_{S}^-$  denotes the Lebesgue measure
on $\G_{S}^-$. Clearly,  $L^2(\G_{S}^-)$ is
dense in $L^2_{w}(\G_{S}^-)$.

 Let
$(\O,\cF,\{\cF_t\}_{t\ge 0},P)$ be a complete
filtered probability space on which a one
dimensional standard Brownian motion
$\{B(t)\}_{t\ge0}$ is defined  such that
$\{\cF_t\}_{t\ge 0}$ is the natural filtration
generated by $\{B(t)\}_{t\ge0}$, augmented by
all the $P$-null sets in $\cF$. Let $H$ be a
Banach space. We denote by $L_{\cF}^2(0,T;H)$
the Banach space consisting of all $H$-valued
$\{\cF_t\}_{t\ge 0}$-adapted processes $X(\cd)$
such that
$\mathbb{E}(|X(\cd)|_{L^2(0,T;H)}^2)<\infty$; by
$L_{\cF}^\infty(0,T;H)$ the Banach space
consisting of all $H$-valued $\{\cF_t\}_{t\ge
0}$-adapted bounded processes;  by
$L_{\cF}^2(\O;C([0,T];H))$ the Banach space
consisting of all $H$-valued $\{\cF_t\}_{t\ge
0}$-adapted continuous processes $X(\cd)$ such
that
$\mathbb{E}(|X(\cd)|_{C([0,T];H)}^2)<\infty$;
and by $C_\cF([0,T];L^2(\O;H))$ the Banach space
consisting of all $H$-valued $\{\cF_t\}_{t\ge
0}$-adapted processes $X(\cd)$ such that
$\mathbb{E}(|X(\cd)|^2)$ is continuous
(similarly, one can define $L^{2}_{\cal
F}(\O;C^{k}([0,T];H))$ and $C^{k}_{\cal
F}([0,T];L^{2}(\O;H))$ for any positive integer
$k$). All of the above spaces are endowed  with
the canonical norm.

The main purpose of this paper is to study the
exact controllability of the following
controlled linear forward stochastic transport
equation:
\begin{equation}\label{csystem1}
\!\!\left\{
\begin{array}{ll}\ds
\2n dy  \! +  \! U\!\cd\! \nabla y dt  \! = \!\[a_1  y +  \!\int_{S^{d-1}}\!\!a_2(t,x,U,V)y(t,x,V)dS^{d-1}(V) +\! f\]dt\\
\ns\ds \hspace{2.7cm} + \big(a_3 y + v\big) dB(t)   &\mbox{ in } (0,T)\t G\t S^{d-1},\\
\ns\ds y  = u  & \mbox{ on } (0,T)\t \G_{S}^-,\\
\ns\ds y(0) = y_0 &\mbox{ in } G\t S^{d-1}.
\end{array}
\right.
\end{equation}
Here and in what follows, $\n$ denotes the
gradient operator with respect to $x$,
$$
\left\{
\begin{array}{ll}\ds
y_0\in L^2(G\times S^{d-1}),\\
\ns\ds a_1\in
L^\infty_{\cF}(0,T;L^\infty(G\times
S^{d-1})),\\
\ns\ds a_2\in
L_\cF^\infty(\O;C([0,T];C(\overline G\times
S^{d-1}\times S^{d-1}))),\\
\ns\ds a_3\in
L^\infty_{\cF}(0,T;L^\infty(G\times
S^{d-1})),\\
\ns\ds f\in L^2_{\cF}(0,T;L^2(G\times
S^{d-1})).
\end{array}
\right.
$$
The boundary control function $u\in
L^2_{\cF}(0,T;L^2_{w}(\G_{S}^-))$ and the
internal control function $v\in
L^2_\cF(0,T;L^2(G\times S^{d-1}))$.

\vspace{0.2cm}

We begin with the definition of solution to the
system \eqref{csystem1}.

\begin{definition}\label{def ft sol}
A solution to  \eqref{csystem1} is a process
$y\in L^2_{\cF}(\O;C([0,T];L^2(G\times
S^{d-1})))$ such that for every $t\in [0,T]$ and
$\phi\in C^1(\overline G\t S^{d-1})$ with
$\phi=0$ on $\G_{S}^+$, it holds that
\begin{equation}\label{def id}
\begin{array}{ll}
\ds \q\int_G\int_{S^{d-1}}
y(t,x,U)\phi(x,U)dS^{d-1}dx -
\int_G\int_{S^{d-1}}
y_0(x,U)\phi(x,U)dS^{d-1}dx\\
\ns\ds \q- \!\int_0^t\!\int_G\!\int_{S^{d-1}}
y(s,x,U)U\cd\nabla\phi(x,U)dS^{d-1}dxds +
\int_0^t\!\int_{\G_{S}^-}\!\! u(s,x,U)\phi(x,U)
U\cd\nu d\G_{S}^- ds
\\ \ns\ds = \int_0^t \int_G\int_{S^{d-1}}
\big[a_1(s,x,U)y(s,x,U) +
f(s,x,U)\big]\phi(x,U)dS^{d-1}dxds \\
\ns\ds \q + \int_0^t\int_G\int_{S^{d-1}} \[\int_{S^{d-1}}a_2(s,x,U,V)y(s,x,V)dS^{d-1}(V)\]\phi(x,U)dS^{d-1}(U)dxds  \\
\ns\ds \q + \int_0^t\int_G\int_{S^{d-1}}
\big[a_3(s,x,U)y(s,x,U) +
v(s,x,U)\big]\phi(x,U)dS^{d-1}dxdB(s),\q
P\mbox{-a.s.}
\end{array}
\end{equation}
\end{definition}

\medskip

In Section 2, we will prove the following
well-posedness result for \eqref{csystem1}.

\medskip

\begin{proposition}\label{well posed}
For each $y_0\in L^2(G\times S^{d-1})$, the
system \eqref{csystem1} admits a unique solution
$y$ such that
\begin{equation}\label{well posed est}
\begin{array}{ll}\ds
\q|y|_{L^2(\O;C([0,T];L^2(G\times S^{d-1})))}\\
\ns\ds \leq e^{Cr_1}\big(|y_0|_{L^2(G\times
S^{d-1})} + |f|_{L^2_{\cF}(0,T;L^2(G\times
S^{d-1}))} +
|u|_{L^2_{\cF}(0,T;L^2_w(\G_{S}^-))} +
|v|_{L^2_{\cF}(0,T;L^2(G\times S^{d-1}))}\big).
\end{array}
\end{equation}
Here  $C>0$ is a constant which is independent
of $y_0$ and
$$
r_1 =
|a_1|^2_{L^\infty_\cF(0,T;L^\infty(G\times
S^{d-1}))} \!+\!
|a_2|_{L_\cF^\infty(\O;C([0,T];C(\overline
G\times S^{d-1}\times S^{d-1})))} \!+\!
|a_3|_{L^\infty_\cF(0,T;L^\infty(G\times
S^{d-1}))}  + 1.
$$
\end{proposition}

\medskip

Now we introduce the notion of exact
controllability for the system \eqref{csystem1}.

\medskip

\begin{definition}\label{exact def}
System \eqref{csystem1} is said to be exactly
controllable at time $T$ if for every  $y_0\in
L^2(G\t S^{d-1})$ and $y_1\in
L^2(\O,\cF_T,P;L^2(G\t S^{d-1}))$, one can find
a pair of controls $(u,v)\in
L^2_\cF(0,T;L^2_w(\G_{S}^-))\t
L^2_{\cF}(0,T;L^2(G\t S^{d-1}))$ such that the
solution $y$ with $y(0)=y_0$ of the system
\eqref{csystem1} satisfies that $y(T) = y_1$.
\end{definition}

\begin{remark}
Since the control $v$ in the diffusion term is
effective in the whole domain, one may expect to
eliminate the randomness of the system
\eqref{csystem1} by taking $v=-a_3y$ and reduce
this system to a controlled random transport
equation. However,  the randomness in
\eqref{csystem1} comes from not only the
stochastic noise $dB$, but also its
coefficients. Although one can take a feedback
control to get rid of the noise term, we still
need to deal with the random coefficients, which
cannot be handled by the classical
controllability theory of deterministic
transport equations.
\end{remark}

We have the following result for the exact
controllability of the system \eqref{csystem1}.

\medskip

\begin{theorem}\label{exact th}
If $T > 2R$, then the system \eqref{csystem1} is
exactly controllable at time $T$.
\end{theorem}

We introduce two controls into the system
\eqref{csystem1}. Moreover, the control $v$ acts
on the whole domain and $T$ needs to be larger
than $2R$. Compared with the deterministic
transport equations, it seems that our choice of
controls is too restrictive. One may consider
the following four weaker cases for designing
the control:

1. Only one control is acted on the system,
that is, $u=0$ or $v=0$ in \eqref{csystem1}.

2. Neither $u$ nor $v$ is zero. But $v=0$ in
$(0,T)\times G_0\times S^{d-1}$, where $G_0$ is
a nonempty open subset of $G$.

3. Two controls are imposed on the system. But
both of them are in the drift term.

4. The time $T<2R$.

It is easy to see that the exact controllability
of \eqref{csystem1} does not hold for the fourth
case. Indeed, if the system \eqref{csystem1}
would be exactly controllable at some time
$T<2R$, then one could deduce the exact
controllability of a deterministic transport
equation on $G$ at time $T$ with a boundary
control acted on $\G_{S}^-$, but this is
obviously impossible. For the other three cases,
according to the controllability result for
deterministic transport equations (see
\cite{Klibanov}), it seems that the
corresponding system should be exactly
controllable. However, as  we shall see later,
it is not  the truth, either.

\medskip
\begin{theorem}\label{non control th}
If $u\equiv 0$ or $v\equiv 0$ in the system
\eqref{csystem1}, then this system is not
exactly controllable at any time $T$.
\end{theorem}
\medskip

Theorem 1.2 indicates that it is necessary to
use two controls  to obtain the desired exact
controllability property for the system
\eqref{csystem1}. Nevertheless,  one may expect
 the exact controllability of \eqref{csystem1}
with the control $v$ (in the diffusion term)
acted only in a proper subdomain of $G$ rather
than the whole domain $G$.  But this is
impossible, either. Indeed, we have the
following negative result.

\medskip

\begin{theorem}\label{non control th1}
 Let $G_0$ be a nonempty open
subset of $G$. If $v\equiv 0$ in $(0,T)\times
G_0\times S^{d-1}$, then the system
\eqref{csystem1} is not exactly controllable at
any time $T$.
\end{theorem}

\medskip

For the third case, we consider the following
controlled equation:
\begin{equation}\label{csystem1.3}
\!\!\!\left\{
\begin{array}{ll}\ds
\2n dy   +   U\cd \nabla y dt   = \[a_1  y\! + \! \int_{S^{d-1}}\!a_2(t,x,U,V)y(t,x,V)dS^{d-1}\!\! +\! f\!+\!\ell\]dt\\
\ns\ds \hspace{3cm} +  a_3 y  dB(t)   &\mbox{ in } (0,T)\t G\t S^{d-1},\\
\ns\ds y  = u  & \mbox{ on } (0,T)\t \G_{S}^-,\\
\ns\ds y(0) = y_0 &\mbox{ in } G\t S^{d-1}.
\end{array}
\right.
\end{equation}
Here $\ell\in L^2_\cF(0,T;L^2(G\times S^{d-1}))$
is another control. Similar to Definition
\ref{exact def}, one can define the exact
controllability of \eqref{csystem1.3}. We have
the following negative result.

\begin{theorem}\label{non control th2}
System \eqref{csystem1.3} is not exactly
controllable for any $T>0$.
\end{theorem}

In order to prove Theorem \ref{exact th}, we
make use of the duality argument.  We obtain the
exact controllability of the system
\eqref{csystem1} by establishing an
observability estimate for  the following
backward stochastic transport equation:
\begin{equation}\label{csystem2}
\left\{ \1n\begin{array}{ll}\ds
\1n \!dz \! + \! U\!\cd\!\nabla zdt \!= \!\!\[b_1 z  + \! \!\int_{S^{d-1}}\!\!\!b_2(t,x,V,U)z(t,x,V)dS^{d-1}(V)\! + \!b_3 Z\] dt \\
\ns\ds \hspace{2.5cm} + (b_4 z +  Z)dB(t)  &\mbox{ in } (0,T)\t G\t S^{d-1},\\
\ns\ds z  = 0 & \mbox{ on } (0,T)\t\G_S^+,\\
\ns\ds z(T) = z_T &\mbox{ in }   G\times
S^{d-1}.
\end{array}
\right.
\end{equation}
Here
$$
\left\{
\begin{array}{ll}\ds
z_T\in L^2(\O,\cF_T,P;L^2(G\times S^{d-1})),\\
\ns\ds b_1\in
L^\infty_{\cF}(0,T;L^\infty(G\times
S^{d-1})),\\
\ns\ds b_2\in
L_\cF^\infty(\O;C([0,T];C(\overline G\times
S^{d-1}\times S^{d-1}))),\\
\ns\ds b_3\in
L^\infty_{\cF}(0,T;L^\infty(G\times
S^{d-1})),\\
\ns\ds b_4\in
L^\infty_{\cF}(0,T;L^\infty(G\times S^{d-1})).
\end{array}
\right.
$$

The definition of  solution to \eqref{csystem2}
is given as follows.

\medskip

\begin{definition}\label{def bt sol}
A solution to the equation \eqref{csystem2} is
a pair of stochastic processes
$$
(z,Z) \in L^2_\cF(\O;C([0,T];\!L^2(G\times
S^{d-1})))\t L^2_{\cF}(0,T;L^2(G\times S^{d-1}))
$$
such that for every $\psi\in C^1(\overline
G\times S^{d-1})$ with $\psi=0$ on $\G_S^-$ and
$t\in [0,T]$, it holds that
\begin{equation}\label{def id1}
\begin{array}{ll}
\ds \q\int_G\int_{S^{d-1}}
z_T(x,U)\psi(x,U)dS^{d-1}dx -
\int_G\int_{S^{d-1}}
z(t,x,U)\psi(x,U)dS^{d-1}dx \\
\ns\ds \q - \int_t^T\int_G\int_{S^{d-1}}
z(s,x,U)U\cd\nabla\psi(x,U)dS^{d-1}dxds
\\ \ns\ds = \int_t^T\int_G\int_{S^{d-1}}
\big[b_1(s,x,U)z(s,x,U) + b_3(s,x,U)Z(s,x,U)\big]\psi(x,U)dS^{d-1}dxds \\
\ns\ds \q + \int_t^T\int_G \int_{S^{d-1}} \[\int_{S^{d-1}}b_2(s,x,V,U)z(s,x,V)dS^{d-1}(V)\]\psi(x,U)dS^{d-1}(U)dxds\\
\ns\ds\q + \int_t^T\int_G\int_{S^{d-1}}
\big[b_4(s,x,U)z(s,x,U) +
Z(s,x,U)\big]\psi(x,U)dS^{d-1}dxdB(s),\;P\mbox{-}\as
\end{array}
\end{equation}
\end{definition}

\medskip

In Section 2, we will establish the following
well-posedness result for \eqref{csystem2}.

\begin{proposition}\label{well posed1}
For any $z_T\in L^2(\O,\cF_T,P;L^2(G\times
S^{d-1}))$, the equation \eqref{csystem2} admits
a unique solution $(z,Z)$ such that
\begin{equation}\label{best1}
|z|_{L^2_\cF(\O;C([0,T];L^2(G\times S^{d-1})))}
+ |Z|_{L^2_\cF(0,T;L^2(G\times S^{d-1}))} \leq
e^{Cr_2} |z_T|_{L^2(\O,\cF_T,P;L^2(G\times
S^{d-1}))},
\end{equation}
where $C$ is a constant which is independent of
$z_T$ and
$$ r_2\=
\sum_{i=1,i\neq
2}^4|b_i|^4_{L^\infty_{\cF}(0,T;L^\infty(G\times
S^{d-1}))} +
|b_2|_{L_\cF^\infty(\O;C([0,T];C(\overline
G\times S^{d-1}\times S^{d-1})))} +1.
$$
\end{proposition}

Now we give the definition of the continuous
observability for the equation \eqref{csystem2}.

\medskip

\begin{definition}\label{exact ob def}
Equation \eqref{csystem2} is said to be
continuously observable in $[0,T]$ if there is a
constant $\cC(b_1,b_2,b_3,b_4)>0$ such that all
solutions of the equation \eqref{csystem2}
satisfy that
\begin{equation}\label{exact ob est}
|z_T|_{L^2(\O,\cF_T,P;L^2(G\times S^{d-1}))}\leq
\cC(b_1,b_2,b_3,b_4)\big(|z|_{L_\cF^2(0,T;L^2_{w}(\G_{S}^-))}
+ |Z|_{L^2_\cF(0,T;L^2(G\times S^{d-1}))}\big).
\end{equation}
\end{definition}

\medskip

The solution $z\in
L^2_{\cF}(\O;C([0,T];L^2(G\times S^{d-1})))$,
hence, it is not obvious that $z|_{\G_S^-}$
belongs to $L^2_\cF(0,T;L^2_{w}(\G_S^-))$. This
is indeed guaranteed by the following regularity
result for \eqref{csystem2}.

\medskip

\begin{proposition}\label{hidden}
Let $(z,Z)$ solve  the equation \eqref{csystem2}
with the terminal state $z_T$. Then
$$
|z|^2_{L^2_\cF(0,T;L^2_{w}(\G_{S}^-))} \leq e^{Cr_2}\mE|z_T|^2_{L^2(G\times
S^{d-1})}.
$$
\end{proposition}

\begin{remark}
The fact that $z|_{\G_{S}^-}\in
L^2_\cF(0,T;L^2_{w}(\G_{S}^-))$ is sometimes
called a hidden regularity property. It does not
follow directly from the classical trace theorem
of Sobolev space.
\end{remark}

It follows from Proposition \ref{hidden}   that
$|z|^2_{L^2_\cF(0,T;L^2_{w}(\G_{S}^-))}$ makes
sense. Now we give the  observability
 result for the equation \eqref{csystem2}.

\medskip

\begin{theorem}\label{th obser}
If $T>2R$, then the equation \eqref{csystem2} is
continuously observable in $[0,T]$.
\end{theorem}

\medskip

In spite of its simple linear form,  the
transport  equation governs many diffusion
processes (see \cite{Landau} for example).
Moreover, it is a linearized Boltzmann equation,
and it is related to the equations of fluid
dynamics such as the Euler and the Navier-Stokes
equations. It is desired to study the stochastic
transport equation  since it is a  model when
the system governed by the transport equation is
perturbed by some stochastic influence.  The
stochastic transport equation is extensively
studied now (see \cite{At1,Cat1,Deck1,Fla1,Pro1}
for example).

The controllability problems for  linear and
nonlinear deterministic transport equations are
well studied in the literature (see
\cite{Coron1,Glass,Leautaud,Klibanov,Rajaram}
and the rich references cited therein).  On the
contrast,  to the author's best knowledge, there
is no published paper addressed to the
controllability of stochastic transport
equations.

Generally speaking, there are three methods to
establish the exact controllability of
deterministic transport equations. The first and
most straightforward one is utilizing the
explicit formula of the solution. By this
method,  for some simple transport equations,
one can explicitly  give a control steering the
system from every given initial state to any
given final state, provided that the time is
large enough. It seems that this method cannot
be used to solve our problem since generally we
do not have the explicit formula for  solutions
to the system \eqref{csystem1}. Nevertheless, we
shall borrow this idea to prove one of our
negative results (i.e., Theorem \ref{non control
th}). The second one is the extension method.
This method was first introduced in
\cite{Russell} to prove the exact
controllability of wave equations. It is
effective to solve the exact controllability
problem for many hyperbolic-type equations.
However, it seems that it is only valid for time
reversible systems. The third and most popular
method is based on the duality between
controllability and observability, via which the
exact controllability problem is reduced to
suitable observability estimate for the dual
system, and the desired observability estimate
is obtained by some global Carleman estimate
(see \cite{Klibanov} for example).

Similar to the deterministic setting, we shall
use a stochastic version of the global Carleman
estimate to derive the inequality \eqref{exact
ob est}. For this, we borrow some idea from the
proof of the observability estimate for the
deterministic transport equations (see
\cite{Klibanov} for example). However, the
stochastic setting will produce some extra
difficulties. We cannot simply mimic the method
in \cite{Klibanov} to solve our problem.

Generally speaking, the nonlocal term, say the
term $\int_{S^{d-1}}
b_2(t,x,V,U)z(t,x,V)dS^{d-1}(V)$ for our
problem, will lead some trouble for obtaining
the observability estimate from the Carleman
estimate, because one cannot simply interchange
the integral operator and the weight function.
However, this will not happen in our case, for
the reason that we choose a weight function
$\theta$ which independent of the variable $U$.

Compared with the extensive results for Carleman
estimate of partial differential equations,
there are a very few works  addressed to its
stochastic counterpart. In \cite{barbu1} and
\cite{Tang-Zhang1}, the authors established some
Carleman type inequalities for forward and
backward stochastic parabolic equations, and via
which the controllability problems for these
equations were addressed. On the other hand, the
authors in \cite{Li-Lu}, \cite{Luqi5} and
\cite{Zhangxu2} obtained some different Carleman
type inequalities for studying  unique
continuation problems for stochastic parabolic
equations. In \cite{Zhangxu3}, a Carleman type
inequality for stochastic wave equations was
first obtained. The result in \cite{Zhangxu3}
was improved in \cite{Luqi3} and \cite{Lu-Zhang}
to solve some inverse problems for stochastic
wave equations. In \cite{Luqi4}, the author got
a Carleman type inequality for stochastic
Schr\"{o}dinger equations and used it to study a
state observation problem for these equations. A
Carleman type inequality for backward stochastic
Schr\"{o}dinger equations was established in
\cite{Luqi1} to prove the exact controllability
of (forward) stochastic Schr\"{o}dinger
equations.

In the literature, in order to obtain the
observability estimate, people usually combine a
Carleman estimate and an Energy estimate (see
\cite{Klibanov} and \cite{Zhangxu1} for
example). In this paper, we deduce the
inequality \eqref{exact ob def} by a new global
Carleman estimate directly (without using the
energy estimate). Indeed, our method even
provide a proof which is simpler than that in
\cite{Klibanov} for the observability estimate
for deterministic transport equations.

The rest of this paper is organized as follows.
In Section 2, we present some preliminary
results, including  the proofs of Proposition
\ref{well posed}-\ref{hidden} and a weighted
identity which is used to prove Theorem \ref{th
obser}. In Section 3, we prove Theorem \ref{th
obser} and in Section 4, we prove Theorem
\ref{exact th}. Finally, Section 5 is addressed
to the proofs of Theorems \ref{non control
th}--\ref{non control th2}.

%%%%%%%%%%%%%%%%%%%%%%%%%%%%%%%%%%%%%%%%%%%%%%%%%%%%%

\section{Some preliminaries}

%%%%%%%%%%%%%%%%%%%%%%%%%%%%%%%%%%%%%%%%%%%%%%%%%%%%%

This section is addressed to present some
preliminary results. We divided it into four
subsections. Proofs of Propositions \ref{well
posed}--\ref{hidden} are given in the first
three subsections. Next, we present a weighted
identity for the stochastic transport operator
$d+U\cdot\nabla dt$, which plays an important
role in establishing the global Carleman
estimate for \eqref{csystem2}.

%%%%%%%%%%%%%%%%%%%%%%%%%%%%%%%%%%%%%%%%%%%%

\subsection{Well-posedness of \eqref{csystem1}}

%%%%%%%%%%%%%%%%%%%%%%%%%%%%%%%%%%%%%%%%%%%%%

In this subsection, we prove Proposition
\ref{well posed}. Equation \eqref{csystem1} is a
nonhomogeneous boundary value problem. Usually,
the well-posedness of such kind of equations is
established in the sense of transposition
solutions (see \cite{Lions1} and \cite{Lions2}
for example). However, fortunately, for our
problem, we can obtain the well-posedness of
\eqref{csystem1} in the context of weak
solution. The key point for doing this is to
establish some suitable a priori estimate (see
the inequality \eqref{well posed ine2} below).

{\it Proof of Proposition \ref{well posed}}\,:
  Let us first
deal with the case in which
\begin{equation}\label{condition1}
\left\{
\begin{array}{ll}\ds y_0 \in L^2(\O,\cF_0,P; H^1(G\t S^{d-1}))
\mbox{ and } y_0 = 0 \mbox{ on } \G_{S}^-, \q
P\mbox{-a.s.},\\
\ns\ds  f,v\in L^2_\cF(0,T;H^1_0(G\t
S^{d-1})),\; u \in Y.
\end{array}
\right.
\end{equation}
Here
$$
\begin{array}{ll}\ds
Y\3n&\ds\=\big\{u:\,u= \tilde u|_{[0,T]\times
\G_S^-} \mbox{ for some } \tilde u\in
 L^2_\cF(\O;C^1([0,T];H^1(G\t S^{d-1}))),\\
\ns&\ds \q\,  \tilde u(0,\cd,\cd)=0 \mbox{ on
}\G_{S}^-,\,\,P\mbox{-a.s.}\big\}.
\end{array}
$$
It is clear that $Y$ is dense in
$L^2_{\cF}(0,T;L^2_w(\G_S^-))$.

Let us consider the following equation:
\begin{equation}\label{csystem1.1}
\left\{
\begin{array}{ll}\ds
\!\!dw\! + U\cd\nabla wdt  \!=\!
\(a_1w \!+\! \int_{S^{d-1}}\!\!a_2(t,x,U,V)w(t,x,V)dS^{d-1}\!\! +\! \tilde f\)dt \\
\ns\ds \hspace{3.0cm} +  (a_3 w + v) dB(t) + a_3 \tilde u dB(t)  &\mbox{ in } \!(0,T)\!\t \!G\t S^{d-1},\\
\ns\ds  w(t,0) = 0 & \mbox{ on }  (0,T)\t\G_{S}^-,\\
\ns\ds  w(0) = y_0 & \mbox{ in }  G\t S^{d-1}.
\end{array}
\right.
\end{equation}
Here
$$
\tilde f = -\tilde u_t-U\cd\nabla \tilde u + a_1
\tilde u + \int_{S^{d-1}}a_2(t,x,U,V)\tilde
u(t,x,V)dS^{d-1} + f.
$$
Clearly, $\tilde f \in L^2_\cF(0,T;H^1(G\times
S^{d-1}))$. Define an unbounded operator $A$ on
$L^2(G\t S^{d-1})$ as follows:
\begin{equation}\label{A}
\left\{
\begin{array}{ll}
\ds D(A) = \left\{  h\in H^1(G\t S^{d-1}):\; h=0 \mbox{ on }\G_{S}^-  \right\},\\
\ns\ds Ah = -U\cd\nabla h, \q\forall\, h \in
D(A).
\end{array}
\right.
\end{equation}
It is an easy matter to see that $D(A)$ is dense
in $L^2(G\t S^{d-1})$ and $A$ is closed.
Furthermore, for every $h\in D(A)$,
$$
(Ah,h)_{L^2(G\t S^{d-1})} = -\int_G\int_{
S^{d-1}} h U\cd\nabla h dS^{d-1}dx =
-\int_{\G_{S}^+}U\cd\nu |h|^2d\G_S^+ \leq 0.
$$
One can easily check that the adjoint
operator of $A$ is
\begin{equation}\label{A star}
\left\{
\begin{array}{ll}
\ds D(A^*) = \left\{  h\in H^1(G\t S^{d-1}):\; h=0 \mbox{ on }\G_{S}^+  \right\},\\
\ns\ds A^*h =  U\cd\nabla h, \q\forall\, h \in
D(A^*).
\end{array}
\right.
\end{equation}
For every $h\in D(A^*)$, it holds that
$$
(A^*h,h)_{L^2\t S^{d-1}} =  \int_G\int_{
S^{d-1}} h U\cd\nabla h dS^{d-1}dx =
\int_{\G_{S}^-}U\cd\nu |h|^2 d\G_{S}^- \leq 0.
$$
Hence, both $A$ and $A^*$ are dissipative
operators. Recalling that $D(A)$ is dense in
$L^2(G\t S^{d-1})$ and $A$ is closed. From the
standard operator semigroup theory (see
\cite[Page 84]{Engel} for example), we conclude
that $A$ generates a $C_0$-semigroup
$\{S(t)\}_{t\geq 0}$ on $L^2(G\t S^{d-1})$ and
$A^*$ generates its dual semigroup
$\{S^*(t)\}_{t\geq 0}$ on $L^2(G\t S^{d-1})$.
Therefore, by the classical theory for
stochastic partial differential equations (see
\cite[Chapter 6]{Prato}),  the system
\eqref{csystem1.1} admits a unique solution
$$w\in L^2_\cF(\O;C([0,T];L^2(G\t
S^{d-1})))\cap L^2_\cF(0,T;D(A))$$ such that
\begin{equation}\label{well posed id1}
\begin{array}{ll}
\ds \q\int_G\int_{S^{d-1}}
w(t,x)\phi(x)dS^{d-1}dx
- \int_G\int_{S^{d-1}}  y_0(x)\phi(x)dS^{d-1}dx\\
\ns\ds \q- \int_0^t\int_G\int_{S^{d-1}}
w(s,x)U\cd\nabla\phi(x)dS^{d-1}dxds
\\ \ns\ds = \int_0^t\int_G\int_{S^{d-1}} \big[a_1(s,x,U) w(s,x,U) + \tilde f(s,x,U)\big]\phi(x,U)dS^{d-1}dxds \\
\ns\ds \q + \int_0^t\int_{S^{d-1}} \[\int_G \int_{S^{d-1}}a_2(s,x,U,V)w(s,x,V)dS^{d-1}(V)\]\phi(x,U)dS^{d-1}(U)dxds \\
\ns\ds \q +
\int_0^t\int_G\int_{S^{d-1}}\!\!\big\{a_3(s,x,U)\big[w(s,x,U)\!+\!\tilde
u(s,x,U)\big]
\!+\! v(s,x,U)\big\}\phi(x,U)dS^{d-1}dxdB(s),\\
\ns\ds\hspace{3cm} P\mbox{-a.s.}, \mbox{ for any
}\phi\in C^1(\overline G\t S^{d-1}) \mbox{ with
} \phi=0 \mbox{ on } \G_{S}^+   \mbox{ and } t
\in [0,T].
\end{array}
\end{equation}

Let
$$
y(t,x,U) = w(t,x,U) + \tilde u(t,x,U), \q \mbox{
for }(t,x,U)\in [0,T]\t G\t S^{d-1}.
$$
Clearly,
$$
y\in L^2_\cF(\O;C([0,T];L^2(G\t S^{d-1})))\cap
L^2_\cF(0,T;H^1(G\times S^{d-1})).
$$
From  \eqref{well posed id1}, we know  that $y$
satisfies
$$
\begin{array}{ll}
\ds \q\int_G\int_{S^{d-1}}
y(t,x,U)\phi(x,U)dS^{d-1}dx -
\int_G\int_{S^{d-1}}
y_0(x,U)\phi(x,U)dS^{d-1}dx \\
\ns\ds \q \!-\!
\!\int_0^t\!\!\int_G\!\int_{S^{d-1}}\!\!\!
y(s,x,U)U\!\cd\!\nabla\phi(x,U)dS^{d-1}\!dxds\!
+ \!\!\int_0^t\!\!\int_G \!\int_{S^{d-1}}\!\!\!\tilde u(s,x,U)U\!\cd\!\nabla\phi(x,U)dS^{d-1}\!dxds \\
\ns\ds =  \!\int_0^t \! \int_G\int_{S^{d-1}}
\big[a_1(s,x,U)y(s,x,U) +  f(s,x,U)- U \cd
\nabla
\tilde u(s,x,U)\big]\phi(x)dS^{d-1}dxds \\
\ns\ds \q + \int_0^t\int_G\int_{S^{d-1}}
\[\int_{S^{d-1}}a_2(s,x,U,V)y(s,x,V)dS^{d-1}(V)\]\phi(x,U)dS^{d-1}(U)dxds  \\
\ns\ds \q + \int_0^t \int_G\int_{S^{d-1}}
\big[a_3(s,x,U)y(s,x,U) +
v(s,x,U)\big]\phi(x,U)dS^{d-1}dxdB(s),\\
\ns\ds\hspace{3cm} P\mbox{-a.s.}, \mbox{ for all
}\phi\in C^1(\overline G\t S^{d-1}) \mbox{ with
} \phi=0 \mbox{ on } \G_{S}^+  \mbox{ and } t\in
[0,T].
\end{array}
$$
Utilizing integration by parts again, we see
that the equality \eqref{def id} holds.
Therefore, $y$ is a solution to the system
\eqref{csystem1} under the assumption
\eqref{condition1}. Furthermore, by means of
It\^{o}'s formula,
\begin{equation}\label{well posed ine1.1}
\!\begin{array}{ll}
\ds \q |y(t)|^2_{L^2(G\times S^{d-1})}\\
\ns\ds =  |y_0|^2_{L^2(G\times S^{d-1})}
 - 2 \int_0^t \int_G \int_{S^{d-1}}  y
U\cd\nabla ydS^{d-1} dxds \\
\ns\ds \q + 2 \int_0^t \int_G \int_{S^{d-1}}
\[\int_{S^{d-1}} a_2 y
dS^{d-1}(V)\]ydS^{d-1}(U)dxds
\\
\ns\ds \q+ \int_0^t\int_G\int_{S^{d-1}}\big[2a_1
y^2 + 2fy + (a_3 y + v)^2 \big] dS^{d-1}dxds\\
\ns\ds \q + 2\int_0^t\int_G\int_{S^{d-1}} y(a_3
y + v) dS^{d-1}dxdB(s).
\end{array}
\end{equation}
This, together with the Burkholder-Davis-Gundy
inequality, implies that
\begin{equation}\label{well posed ine1}
\!\begin{array}{ll}
\ds \q\mE \sup_{s\in [0,t]}|y(s)|^2_{L^2(G\times S^{d-1})}\\
\ns\ds   \leq  |y_0|^2_{L^2(G\times S^{d-1})}
\!\!-\! \mE\!\int_0^t\!\int_{\G_{S}^-}\!U\cd\nu
u^2 d\G_S^- ds \!+\!
2\mE\!\int_0^t\!\int_G\!|a_2|_{C(\overline
G\times S^{d-1}\times
S^{d-1})}\!\int_{S^{d-1}}\!\! y^2dS^{d-1}dxds
\\
\ns\ds \q + 4\mE\int_0^t\int_G\int_{S^{d-1}}
\big[a_1y^2 + y^2 + f^2 + a_3^2 y^2 + v^2\big]
dS^{d-1}dxds \\
\ns\ds \leq  |y_0|^2_{L^2(G\times S^{d-1})}\!
+4r_1\mE\!\int_0^t\!  \big[\sup_{\si\in [0,s]}
|y(\si)|^2_{L^2(G\times S^{d-1})}\big] ds -
\mE\int_0^t\!\int_{\G_{S}^-}\! U\cd\nu u^2
d\G_{S}^-
ds \\
\ns\ds \q +
4\mE\!\int_0^t\!\int_G\!\int_{S^{d-1}}\!\big(
f^2 \!+\! v^2\big)dS^{d-1}dxds.
\end{array}
\end{equation}
Hence, by Gronwall's inequality, we obtain that
\begin{equation}\label{well posed ine2}
\begin{array}{ll}\ds \q
|y|_{L^2_\cF(\O;C([0,T];L^2(G\times S^{d-1})))}\\
\ns\ds \leq e^{Cr_1}\big( |y_0|_{L^2(G\times
S^{d-1})}\! +\!
|u|_{L^2_{\cF}(0,T;L^2_{w}(\G_{S}^-))} \!+\!
|f|_{L^2_{\cF}(0,T;L^2(G\times S^{d-1}))} \!+\!
|v|_{L^2_{\cF}(0,T;L^2(G\times S^{d-1}))} \big).
\end{array}
\end{equation}

By a similar argument, we can show that
if
$$
\begin{array}{ll}\ds
(\hat y_0,\hat u,\hat f,\hat v)\3n&\ds\in D(A)
\t Y \t  L^2_\cF(0,T;H_0^1(G\times S^{d-1}))\t
L^2_\cF(0,T;H_0^1(G\times S^{d-1}))
\end{array}
$$
and
$$
\begin{array}{ll}\ds
(\bar y_0, \bar u,\bar f,\bar v)\3n&\ds\in D(A)
\t Y\t  L^2_\cF(0,T;H_0^1(G\times S^{d-1}))\t
L^2_\cF(0,T;H_0^1(G\times S^{d-1})),
\end{array}
$$
then we can find corresponding solutions
$$\hat y,\bar y\in
L^2_\cF(\O;C([0,T];L^2(G\times S^{d-1}))) \cap
L^2_\cF(0,T;H^1(G\times S^{d-1}))$$ such that
$$
\begin{array}{ll}\ds
\q|\hat y -\bar
y|_{L^2_\cF(\O;C([0,T];L^2(G\times S^{d-1})))}\\
\ns\ds \leq e^{Cr_1}\big( |\hat y_0-\bar
y_0|_{L^2(G\times S^{d-1})} + |\hat u-\bar
u|_{L^2_{\cF}(0,T;L^2_{w}(\G_{S}^-))} + |\hat f
- \bar f|_{L^2_\cF(0,T;L^2(G\times S^{d-1}))}
\\
\ns\ds \q + |\hat v- \bar
v|_{L^2_\cF(0,T;L^2(G\times S^{d-1}))} \big).
\end{array}
$$

Now for $y_0 \in L^2(G\times S^{d-1})$, $u\in
L^2_{\cF}(0,T;L^2_{w}(\G_{S}^-))$, $f, v\in
L^2_\cF(0,T;L^2(G\times S^{d-1}))$, let us
choose
$$
\begin{array}{ll}\ds
\{y_0^n\}_{n=1}^{+\infty}\subset  D(A),
\q\{u^n\}_{n=1}^{+\infty}\subset Y,\q
\{f^n\}_{n=1}^{+\infty}\subset
L^2_\cF(0,T;H_0^1(G\times
S^{d-1})),\\
\ns\ds \{v^n\}_{n=1}^{+\infty}\subset
L^2_\cF(0,T;$ $H_0^1(G\times S^{d-1})),
\end{array}
$$
such that
\begin{equation}\label{y0n}
\left\{
\begin{array}{ll}\ds
\lim_{n\to\infty} y_0^n = y_0 \mbox{ in }
L^2(G\times S^{d-1});\\
\ns\ds \lim_{n\to\infty} u^n = u \mbox{ in }
L^2_{\cF}(0,T;L^2_{w}(\G_{S}^-));\\
\ns\ds \lim_{n\to\infty} f^n = f \mbox{ in }
L^2_{\cF}(0,T;L^2(G\times S^{d-1})); \\
\ns\ds \lim_{n\to\infty} v^n = v \mbox{ in }
L^2_{\cF}(0,T;L^2(G\times S^{d-1})).
\end{array}
\right.
\end{equation}
For every given $(y_0^n, u^n,f^n,v^n)$, by the
argument above, we know that there is a unique
solution $y_n(\cdot,\cdot)$ to the system
\eqref{csystem1}, which satisfies
\begin{equation}\label{well posed id3}
\begin{array}{ll}
\ds \q\int_G\int_{S^{d-1}}
y_n(t,x,U)\phi(x)dS^{d-1}dx -
\int_G\int_{S^{d-1}}
y^n_0(x,U)\phi(x,U)dS^{d-1}dx \\
\ns\ds \q -\!
\int_0^t\!\int_G\!\int_{S^{d-1}}\!\!\!
y_n(s,x,U)U\!\cd\!\nabla\phi(x,U)dS^{d-1}dxds
\!-\! \int_0^\tau\!\!\int_{\G_{S}^-}\!
U\!\cd\!\nu u^n(s,x,U)\phi(x,U)d\G_{S}^-ds \\
\ns\ds = \int_0^t \int_G\int_{S^{d-1}}
\big[a_1(s,x,U)y_n(s,x,U) +
f^n(s,x,U) \big]\phi(x)dS^{d-1}dxds \\
\ns\ds \q + \int_0^t\int_G\int_{S^{d-1}}
\[\int_{S^{d-1}}a_2(s,x,U,V) y_n(s,x,V)dS^{d-1}(V)\]\phi(x,U)dS^{d-1}(U)dxds  \\
\ns\ds \q + \int_0^t \int_G\int_{S^{d-1}}
\big[a_3(s,x,U)y_n(s,x,U) +
v^n(s,x,U)\big]\phi(x,U)dS^{d-1}dxdB(s),\\
\ns\ds\hspace{3cm} P\mbox{-a.s.}, \mbox{ for any
}\phi\in C^1(\overline G\t S^{d-1}) \mbox{ with
} \phi=0 \mbox{ on } \G_{S}^+ \mbox{ and } \tau
\in [0,T],
\end{array}
\end{equation}
and
\begin{equation}\label{well posed ine2.1}
\begin{array}{ll}\ds
\q
|y_n|_{L^2_\cF(\O;C([0,T];L^2(G\times S^{d-1})))}\\
\ns\ds \leq \! e^{Cr_1}\!\big(|y_0^n|_{L^2(G\t
S^{d-1})}\!  +\!
|u^n|_{L^2_{\cF}(0,T;L^2_{w}(\G_{S}^-))} \! +\!
|f^n|_{L^2_{\cF}(0,T;L^2(G\times S^{d-1}))} \!
+\! |v^n|_{L^2_{\cF}(0,T;L^2(G\times S^{d-1}))}
\big).
\end{array}
\end{equation}
Further, for any $m,n\in \dbN$, we have
\begin{equation}\label{well posed ine4}
\begin{array}{ll}\ds
\q|y_n -
y_m|_{L^2_\cF(\O;C([0,T];L^2(G\times S^{d-1})))}\\
\ns\ds \leq e^{Cr_1}\big( | y_0^n-
y_0^m|_{L^2(G\times S^{d-1})} + |u^n-
u^m|_{L^2_{\cF}(0,T;L^2_{w}(\G_{S}^-))} + | f^n
- f^m|_{L^2_\cF(0,T;L^2(G\times S^{d-1}))}
\\
\ns\ds \q + |v^n- v^m|_{L^2_\cF(0,T;L^2(G\times
S^{d-1}))} \big).
\end{array}
\end{equation}
From \eqref{y0n} and \eqref{well posed ine4}, we
obtain that $\{y_n\}_{n=1}^{+\infty}$ is a
Cauchy sequence in
$L^2_\cF(\O;C([0,T];L^2(G\times S^{d-1})))$.
Hence, there exists a unique $y\in
L^2_\cF(\O;C([0,T];L^2(G\times S^{d-1})))$ such
that
\begin{equation}\label{well posed4}
y_n \to y \mbox{ in }
L^2_\cF(\O;C([0,T];L^2(G\times S^{d-1})))\;
\mbox{ as } n\to +\infty.
\end{equation}
Combining \eqref{well posed id3} and \eqref{well
posed4}, we find that $y$ satisfies \eqref{def
id}. Hence, $y$ is a solution to the system
\eqref{csystem1}.

Further, from \eqref{well posed ine2.1} and
\eqref{well posed4}, we obtain that $y$
satisfies the inequality \eqref{well posed est}.

The uniqueness of the solution to
\eqref{csystem1} follows from \eqref{well posed
est} immediately. This completes the proof of
Proposition \ref{well posed}.
\endpf

%%%%%%%%%%%%%%%%%%%%%%%%%%%%%%%%%%%%%%%%%

\subsection{Well-posedness of \eqref{csystem2}}

%%%%%%%%%%%%%%%%%%%%%%%%%%%%%%%%%%%%%%%%%

This subsection is devoted to a proof of
Proposition \ref{well posed1}.

We first recall the definition of the mild
solution to backward stochastic evolution
equations.

Let $X$ be a Hilbert space and
$\cA:D(\cA)\subset X \to X$ be a linear operator
which generates a $C_0$-semigroup
$\{\cS(t)\}_{t\geq 0}$ on $X$. Let
$F_1:[0,T]\times X\times X \to X$ satisfy that
\begin{itemize}
  \item there exists an $L_1>0$ such that
  $$
\!\!\!\!\!\!\!\!\!\!|F_1(t,\eta_1,\eta_2) -
F_1(t,\hat \eta_1,\widehat \eta_2)|_X \leq
L_1(|\eta_1-\hat \eta_1|_X + |\eta_2-\widehat
\eta_2|_X) \mbox{ for all } t\!\in [0,T],\;
\eta_1,\hat \eta_1, \eta_2, \widehat \eta_2 \in
X;
 $$
  \item $F_1(\cd,0,0)\in L^2(0,T;X)$.
\end{itemize}
Let $F_2(\cd,\cd):[0,T]\times X \to X$ satisfy
that
\begin{itemize}
 \item there exists an $L_2>0$ such that
 $$
|F_2(t, \eta_1) - F_2(t,\hat \eta_1)|_X \leq
L_2|\phi-\hat \phi|_X\mbox{ for all } t\in
[0,T],\; \eta_1,\hat \eta_1 \in X;
  $$
  \item $F_2(\cd,0)\in L^2(0,T;X)$.
\end{itemize}

Consider the following backward stochastic
evolution equation
\begin{equation}\label{bsystem1}
\left\{
\begin{array}{ll}\ds
d\phi = -\big[\cA\phi(t) +
F_1(t,\phi(t),\Phi(t))\big]dt -
\big[F_2(t,\phi(t))+\Phi(t)\big]dB(t)
&\mbox{ in } [0,T],\\
\ns\ds \phi(T)=\phi_T,
\end{array}
\right.
\end{equation}
where $\phi_T\in L^2(\O,\cF_T,P;X)$.

A pair of processes $(\phi,\Phi)\in
L^2_\cF(\O;C([0,T];X))\times L^2_\cF(0,T;X)$ is
a mild solution of \eqref{bsystem1} if for all
$t\in [0,T]$, they satisfy that
\begin{equation}\label{mild}
\begin{array}{ll}\ds
\phi(t)\3n&\ds =\cS(T-t)\phi_T + \int_t^T
\cS(s-t)F_1(s,\phi(s),\Phi(s))ds\\
\ns&\ds \q + \int_t^T
\cS(s-t)\big[F_2(s,\phi(s))+\Phi(s)\big]dB(s),\q
P\mbox{-}\as
\end{array}
\end{equation}
\begin{lemma}\cite[Theorem 9]{Mah1}\label{lm2}
The equation \eqref{bsystem1} admits a unique
mild solution $(\phi,\Phi)$.
\end{lemma}

We are now in a position to prove Proposition
\ref{well posed1}.

\vspace{0.15cm}

{\it Proof of Proposition \ref{well
posed1}}\,:\, Let $X=L^2(G\times S^{d-1})$,
$\cA=A^*$,
$$
\left\{
\begin{array}{ll}\ds
 F_1(t,\phi,\Phi) =  -\[b_1 \phi  +
\!\int_{S^{d-1}}\!
b_2(t,x,V,U)\phi(t,x,V)dS^{d-1}(V)
+  b_3 \Phi\],\\
\ns\ds F_2(t,\phi) = -b_4\phi.
\end{array}
\right.
$$
We have $\cS(t)=S^*(t)$. By Lemma \ref{lm2}, we
conclude that \eqref{csystem2} admits a unique
mild solution $(z,Z)$ such that
\begin{equation}\label{10.15-eq1}
\begin{array}{ll}\ds
z(t)\3n&\ds=S^*(T-t)z_T - \int_t^T S^*(s-t)\[b_1
z + \!  \int_{S^{d-1}}\!
b_2(s,x,V,U)z(s,x,V)dS^{d-1}(V)\!
+ \!b_3 Z\]ds\\
\ns&\ds \q - \int_t^T S^*(s-t)\big(b_4z +
Z\big)dB(s),\q P\mbox{-}\as
\end{array}
\end{equation}
From \eqref{10.15-eq1}, for any $\psi\in
C^1(\overline G\times S^{d-1})$ with $\psi=0$ on
$\G_S^-$, we have that
\begin{equation}\label{10.15-eq2}
\begin{array}{ll}\ds
\q\langle z(t),A\psi\rangle_{L^2(G\times
S^{d-1})}
\\
\ns\ds = \langle
S^*(T\!-t)z_T,A\psi\rangle_{L^2(G\times
S^{d-1})}\! + \!\int_t^T \!\!\Big\langle
S^*(s\!-t) F_1(s,z(s),Z(s)), A\psi
\Big\rangle_{L^2(G\times S^{d-1})}ds\\
\ns\ds\q - \int_t^T \big\langle
S^*(s-t)\big(b_4z +
Z\big),A\psi\big\rangle_{L^2(G\times S^{d-1})}
dB(s)\\
\ns\ds = \langle
z_T,S(T-t)A\psi\rangle_{L^2(G\times S^{d-1})} +
\int_t^T \Big\langle F_1(s,z(s),Z(s)),
S(s-t)A\psi \Big\rangle_{L^2(G\times S^{d-1})}ds
\\
\ns\ds \q- \int_t^T \big\langle b_4z +
Z,S(s-t)A\psi\big\rangle_{L^2(G\times S^{d-1})}
dB(s)\\
\ns\ds \= I_1 + I_2 - I_3.
\end{array}
\end{equation}
Integrating \eqref{10.15-eq2} from $t$ to $T$,
we obtain that
\begin{equation}\label{10.15-eq3}
\int_t^T \langle z(t),A\psi\rangle_{L^2(G\times
S^{d-1})}ds = \int_t^T (I_1 + I_2 - I_3)ds.
\end{equation}
Clearly,
\begin{equation}\label{10.15-eq4}
\begin{array}{ll}\ds
\int_t^T I_1 ds \3n&\ds = \int_t^T \langle
z_T,S(T-s)A\psi\rangle_{L^2(G\times
S^{d-1})}ds\\
\ns&\ds = \langle
z_T,S(T-t)\psi\rangle_{L^2(G\times S^{d-1})} -
\langle z_T, \psi\rangle_{L^2(G\times S^{d-1})}.
\end{array}
\end{equation}
By the  Fubini's theorem, we have that
\begin{equation}\label{10.15-eq5}
\begin{array}{ll}\ds
\int_t^T I_2 ds \3n&\ds= \int_t^T \int_s^T
\Big\langle F_1(r,z(r),Z(r)), S(r-s)A\psi
\Big\rangle_{L^2(G\times S^{d-1})}drds\\
\ns&\ds = \int_t^T \Big\langle F_1(r,z(r),Z(r)),
\int_t^r S(r-s)A\psi ds
\Big\rangle_{L^2(G\times S^{d-1})}dr\\
\ns&\ds = \int_t^T \Big\langle F_1(r,z(r),Z(r)),
S(r-t) \psi - \psi \Big\rangle_{L^2(G\times
S^{d-1})}dr \\
\ns&\ds = \int_t^T \Big\langle S^*(r-t)
F_1(r,z(r),Z(r)), \psi  \Big\rangle_{L^2(G\times
S^{d-1})}dr\\
\ns&\ds \q - \int_t^T \Big\langle
F_1(r,z(r),Z(r)), \psi \Big\rangle_{L^2(G\times
S^{d-1})}dr,
\end{array}
\end{equation}
and by the stochastic Fubini's theorem (see
\cite[page 109]{Prato} for example), we find
that
\begin{equation}\label{10.15-eq6}
\begin{array}{ll}\ds
\int_t^T I_3 ds \3n&\ds= \int_t^T \int_s^T
\big\langle b_4(r)z(r) +
Z(r),S(r-s)A\psi\big\rangle_{L^2(G\times
S^{d-1})} drdB(s)\\
\ns&\ds = \int_t^T \Big\langle b_4(r)z(r) +
Z(r), \int_t^r S(r-s)A\psi ds
\Big\rangle_{L^2(G\times S^{d-1})}dB(r)\\
\ns&\ds = \int_t^T \Big\langle b_4(r)z(r) +
Z(r),  S(r-t) \psi -\psi
\Big\rangle_{L^2(G\times S^{d-1})}dB(r)\\
\ns&\ds = \int_t^T \Big\langle S^*(r-t)
\big[b_4(r)z(r) + Z(r)\big],  \psi
\Big\rangle_{L^2(G\times S^{d-1})}dB(r)\\
\ns&\ds \q - \int_t^T \Big\langle b_4(r)z(r) +
Z(r), \psi \Big\rangle_{L^2(G\times
S^{d-1})}dB(r).
\end{array}
\end{equation}
From \eqref{10.15-eq2}--\eqref{10.15-eq6}, we
obtain that $(z,Z)$ satisfies \eqref{def id1}.

The proof of the inequality \eqref{best1} is
very similar to the one of \eqref{well posed
est}. Indeed, by It\^{o}'s formula, we can
easily obtain that
\begin{equation}\label{10.15-eq7}
\begin{array}{ll}\ds
\q|z_T|_{L^2(G\times S^{d-1})}^2-|z(t)|_{L^2(G\times S^{d-1})}^2 \\
\ns\ds \geq 2 \int_t^T\int_G \int_{S^{d-1}}
z\[b_1 z + \! \int_{S^{d-1}}\!\!
b_2(r,x,V,U)z(t,x,V)dS^{d-1}(V) +  b_3 Z\]dS^{d-1}dxds \\
\ns\ds \q + 2\int_t^T\int_G \int_{S^{d-1}}
z\big(b_4 z +  Z\big)dS^{d-1}dxdB(s) +
\int_t^T\int_G \int_{S^{d-1}}\big(b_4 z +
Z\big)^2 dS^{d-1}dxds.
\end{array}
\end{equation}
By Burkholder-Davis-Gundy inequality, we find
that
\begin{equation}\label{10.15-eq8}
\begin{array}{ll}\ds
\q \mE\sup_{s\in [t,T]}|z(t)|_{L^2(G\times
S^{d-1})}^2 + |Z|^2_{L^2_\cF(t,T;L^2(G\times
S^{d-1}))}
\\
\ns\ds \leq |z_T|_{L^2(G\times S^{d-1})}^2+Cr_2
\mE\int_t^T |z(s)|^2_{L^2(G\times S^{d-1})}.
\end{array}
\end{equation}
Then, by the Gronwall's inequality, we get
\eqref{best1} immediately. The uniqueness of the
solution follows from the inequality
\eqref{best1}. This completes the proof of
Proposition \ref{well posed1}.
\endpf

%%%%%%%%%%%%%%%%%%%%%%%%%%%%%%%%%%%%%%%%%%%%%%%

\subsection{Hidden regularity for solutions to backward stochastic transport equations}

%%%%%%%%%%%%%%%%%%%%%%%%%%%%%%%%%%%%%%%%%%%%%%%

In this subsection, we   give a proof of
Proposition \ref{hidden}.

\vspace{0.2cm}

{\it Proof of Proposition \ref{hidden}}\,: The
proof  is almost standard. Here  we give it for
the sake of completeness. Let
$$
\cX \= \big\{h\in H^1(G\times S^{d-1}):\, h=0
\mbox{ on } \G_{S}^+ \big\}.
$$
Following the proof of Proposition \ref{well
posed1} (for this, one needs numerous but small
changes), one can show that if $z_T\in
L^2(\O,\cF_T,P;\cX)$, then the solution
$$
(z,Z)\in \big(L^2_\cF(\O;C([0,T];L^2(G\times
S^{d-1})))\cap L^2_\cF(0,T;\cX)\big)\t
L^2_\cF(0,T;L^2(G\times S^{d-1})).
$$
Then, by It\^{o}'s formula, we see that
\begin{equation}\label{hidden1}
\begin{array}{ll}\ds
\q\mE|z_T|^2_{L^2(G\times S^{d-1})} - |z(0)|^2_{L^2(G\times S^{d-1})} \\
\ns\ds = \!-\mE\!\int_0^T\!\int_G
\int_{S^{d-1}}\!\! zU\!\cd\!\nabla z
dS^{d-1}dxdt \!+\! \mE\!\int_0^T\!\int_G
\int_{S^{d-1}}\!\! \big[ 2z(b_1 z
\!+\! b_3 Z) \!+\! (b_4 z \!+\! Z)^2 \big]dS^{d-1}dxdt\\
\ns\ds \q + \mE\int_0^T\int_G\int_{S^{d-1}}
z(t,x,U)\[\int_{S^{d-1}}
b_2(t,x,U,V)z(t,x,V)dS^{d-1}(V)
\]dS^{d-1}(U)dxdt.
\end{array}
\end{equation}
Therefore, we find that
\begin{equation}\label{hidden2}
\begin{array}{ll}\ds
\q-\mE\int_0^T\int_{\G_{S}^-} U\cd\nu z^2d\G_{S}^- dt \\
\ns\ds  = \mE|z_T|^2_{L^2(G\times S^{d-1})}
\!\!- \!  |z(0)|^2_{L^2(G\times S^{d-1})}\!
\!+\! \mE\int_0^T \!\! \!\int_G \int_{S^{d-1}}
\!\!\!\big[ 2z(b_1 z
\!+\! b_3 Z) \!+\! (b_4 z \!+\! Z)^2 \big]dS^{d-1}dxdt\\
\ns\ds \q + \mE\int_0^T\int_G\int_{S^{d-1}}
z(t,x,U)\[\int_{S^{d-1}}
b_2(t,x,U,V)z(t,x,V)dS^{d-1}(V)
\]dS^{d-1}(U)dxdt\\
\ns\ds \leq e^{Cr_2}\mE|z_T|^2_{L^2(G\times
S^{d-1})}.
\end{array}
\end{equation}
For any $z_T \in L^2(\O,\cF_T,P;L^2(G\times
S^{d-1}))$, we can find a sequence
$\{z_T^{(n)}\}_{n=1}^\infty\subset
L^2(\O,\cF_T,P;\cX)$ such that
$$
\lim_{n\to\infty}z_T^{(n)} = z_T \mbox{ in
}L^2(\O,\cF_T,P;L^2(G\times S^{d-1})).
$$
Hence, we know that the inequality
\eqref{hidden2} also holds for $z_T \in
L^2(\O,\cF_T,P;L^2(G\times S^{d-1}))$. \endpf

%%%%%%%%%%%%%%%%%%%%%%%%%%%%%%%%%%%%%%%%%%%%%%%%%

\subsection{Identity for a stochastic transport operator}

%%%%%%%%%%%%%%%%%%%%%%%%%%%%%%%%%%%%%%%%%%%%%%%%%

In this subsection,  we introduce a weighted
identity for the stochastic transport operator
$d+U\cdot\nabla dt$, which will play a key role
in the proof of Theorem \ref{th obser}.  Let
$\l>0$, and let $0<c<1$ such that $cT>2R$. Put
\begin{equation}\label{weight1} l = \l \[
|x|^2  - c \(t-\frac{T}{2}\)^2 \] \;\mbox{ and }
\;\theta = e^{l}.
\end{equation}

We have the following weighted identity involving
$\theta$ and $l$.

\medskip

\begin{proposition}\label{prop identity}
Assume that $q$ is an $H^1(\dbR^n)\times
L^2(S^{d-1})$-valued continuous semi-martingale.
Put $p=\theta q$. We have the following equality
\begin{equation}\label{identity}
\begin{array}{lll}
\ds \q -\theta(l_t + U\cd\nabla l)p\big[ dq + U\cd\nabla q dt \big]\\
\ns\ds = -\frac{1}{2}d\big[ (l_t + U\cd\nabla
l)p^2 \big]
- \frac{1}{2} U\cd\nabla\big[(l_t + U\cd\nabla l)p^2\big] + \frac{1}{2}\big[l_{tt} + U\cd\nabla (U\cd\nabla l)\\
\ns\ds \q +  2U\cd\nabla l_{t}  \big] p^2 +
\frac{1}{2}(l_t + U\cd\nabla l)(dp)^2 + (l_t +
U\cd\nabla l)^2 p^2.
\end{array}
\end{equation}
\end{proposition}

\medskip

{\it Proof of Proposition \ref{prop
identity}}\,:  By the definition of $p$, we
have
$$
\theta(dq+U\cd\nabla q) =\theta d(\theta^{-1} p)
+ \theta U\cd\nabla(\theta^{-1} p) = dp +
U\cd\nabla p - (l_t + U\cd\nabla l)p.
$$
Thus,
\begin{equation}\label{identity.1}
\begin{array}{ll}
\q -\theta(l_t + U\cd\nabla l)p\big( dq + U\cd\nabla q  \big)\\
\ns\ds = -(l_t +U\cd\nabla l )p \big[ dp + U\cd\nabla p - (l_t + U\cd\nabla l)p \big]\\
\ns\ds = -(l_t + U\cd\nabla l)p (dp +
U\cd\nabla p) + (l_t + U\cd\nabla l)^2p^2.
\end{array}
\end{equation}
It is easy to see that
\begin{equation}\label{identity.2}
\left\{
\begin{array}{lll}\ds
-l_t pdp = -\frac{1}{2}d(l_t p^2) + \frac{1}{2}l_{tt}p^2 + \frac{1}{2}l_t (dp)^2,\\
\ns\ds -U\cd\nabla l pdp = -\frac{1}{2}d(U\cd\nabla l p^2) + \frac{1}{2}(U\cd\nabla l)_t p^2 + \frac{1}{2}U\cd\nabla l (dp)^2,\\
\ns\ds -l_{t} p U\cd\nabla p = -\frac{1}{2}U\cd\nabla(l_t p^2) + \frac{1}{2}U\cd\nabla l_{t}p^2,\\
\ns\ds -U\cd\nabla l pU\cd\nabla p=
-\frac{1}{2}U\cd\nabla(U\cd\nabla l p^2)  +
\frac{1}{2}U\cd\nabla(U\cd\nabla l)p^2.
\end{array}
\right.
\end{equation}
From  \eqref{identity.1} and
\eqref{identity.2}, we obtain the equality
\eqref{identity}.
\endpf

%%%%%%%%%%%%%%%%%%%%%%%%%%%%%%%%%%%%%%%%%%%%%%

\section{Proof of Theorem \ref{th obser}}

%%%%%%%%%%%%%%%%%%%%%%%%%%%%%%%%%%%%%%%%

This section is devoted to proving Theorem
\ref{th obser} by means of  a suitable global
Carleman estimate for the equation
\eqref{csystem2}.

\vspace{0.2cm}

{\it Proof of Theorem \ref{th obser}}\,: To
begin with, applying Proposition \ref{prop
identity} to the equation \eqref{csystem2} with
$v=z$, integrating \eqref{identity} on $(0,T)\t
G\t S^{d-1}$ and using integration by parts, and
taking expectation, we get that
\begin{equation}\label{bobeq1}
\begin{array}{lll}\ds
\q-2\mE \int_0^T\int_G\int_{S^{d-1}} \theta^2(l_t + U\cd\nabla l)z (dz + U\cd\nabla z dt)dS^{d-1}dxdt\\
\ns\ds =\l \mE\int_G\!\int_{S^{d-1}}\!\!
(cT\!-\!2U\cd x) \theta^2(T) z^2(T)dS^{d-1}dx
+\l \!\int_G\!\int_{S^{d-1}}\!\! (cT \!+\!
2U\!\cd\! x)\theta^2(0)
z^2(0)dS^{d-1}dx\\
\ns\ds \q + \l\mE\int_0^T \!\int_{\G_{S}^-}
U\cd\nu \big[c(T\!-\!2t)\!-\!2 U\!\cd\!
x\big]\theta^2 z^2 d\G_{S}^- dt \!+\!
2(1\!-\!c)\l\mE\!\int_0^T\!\int_G\!\int_{S^{d-1}}
\theta^2 z^2 dS^{d-1}dxdt
\\
\ns\ds \q +
\mE\!\int_0^T\!\!\!\!\int_G\!\int_{S^{d-1}}
\!\!\!\theta^2 (l_t \!+\! U\!\cd\!\nabla l)
(b_4z\!+\!Z)^2dS^{d-1}dxdt \!+\!
2\mE\!\int_0^T\!\!\!\!\int_G\!
\int_{S^{d-1}}\!\!\!\theta^2 (l_t \!+\!
U\!\cd\!\nabla l)^2 z^2dS^{d-1}dxdt.
\end{array}
\end{equation}
By virtue of that $z$ solves the equation
\eqref{csystem2}, we see that
\begin{equation}\label{bobeq2}
\begin{array}{lll}\ds
\q-2\mE \int_0^T\int_G\int_{S^{d-1}} \theta^2(l_t + U\cd\nabla l) z (dz + U\cd\nabla zdt)dS^{d-1}dxdt\\
\ns\ds =  \!2\mE\!
\!\int_0^T\!\!\!\int_G\!\int_{S^{d-1}}\!\!\!\!
\theta^2(l_t \!+\! U\!\cd\!\nabla l) z \(\! b_1z
\!+\!\!
\int_{S^{d-1}}\!\!\!\!b_2(t,x,\!U,\!V)z(t,x,\!V)dS^{d-1}(V)\!
+\! b_3Z\!\) dS^{d-1}(U)dxdt
\\
\ns\ds\leq \mE\int_0^T\!\!\int_G\!
\int_{S^{d-1}}\!\theta^2 (l_t + U\cd\nabla l)^2
z^2dS^{d-1}dxdt +
3\mE\int_0^T\!\!\int_G\!\int_{S^{d-1}}\!
\theta^2 (b_1^2 z^2 + b_3^2 Z^2) dS^{d-1}dxdt \\
\ns\ds \q + 3\mE\int_0^T\int_G\int_{S^{d-1}}
\theta^2\|\int_{S^{d-1}}b_2(t,x,U,V)z(t,x,V)dS^{d-1}(V)\|^2dS^{d-1}(U)dxdt.
\end{array}
\end{equation}
This, together with the equality
\eqref{bobeq1}, implies that
\begin{equation}\label{bobeq3}
\begin{array}{lll}\ds
\q \l \mE\!\int_G\!\int_{S^{d-1}}\!\!\!
(cT\!\!-\!2U\cd x) \theta^2(T) z^2(T)dS^{d-1}dx
\!+\!\l \!\int_G\!\int_{S^{d-1}}\!\!\! (cT\! +\! 2U\cd x)\theta^2(0) z^2(0)dS^{d-1}dx\\
\ns\ds\q  + 2(\!1 \!-\! c)\l\mE\!\int_0^T\!\!
\int_G\!\int_{S^{d-1}}\!\!\! \theta^2 z^2
dS^{d-1}dxdt \!+\! \mE\int_0^T \!\!\!\int_G\!
\int_{S^{d-1}}\!\!\theta^2 (l_t \!+\!
U\!\cd\!\nabla l) (b_4 z \!+\! Z)^2dS^{d-1}dxdt
\\
\ns\ds \q +\mE\int_0^T \!\int_G\int_{S^{d-1}}
\theta^2 (l_t  +  U\cd\nabla l)^2
z^2dS^{d-1}dxdt
\\
\ns\ds \leq \!3\mE \!\int_0^T\!
\!\int_G\int_{S^{d-1}} \!\!\theta^2 (b_1^2 z^2
\!+\! b_3^2 Z^2)dS^{d-1}dxdt \!-\! \l
\mE\int_0^T\!\! \int_{\G_{S}^-}\! U\cd\nu
\big[c(T\!-\!2t)\!-\!2 U\cd x\big]\theta^2 z^2
d\G_{S}^-
dt\\
\ns\ds \q + 3
|b_2|^2_{L^\infty_\cF(\O;L^\infty(0,T;C(\overline
G\times S^{d-1}\times S^{d-1})))}\mE\int_0^T
\int_G\int_{S^{d-1}} \theta^2 z^2dS^{d-1}dxdt.
\end{array}
\end{equation}
Since
$$
\begin{array}{ll}\ds
\q \mE\int_0^T \int_G \int_{S^{d-1}}\theta^2 (l_t + U\cd\nabla l) (b_4 z + Z)^2dS^{d-1}dxdt
\\
\ns\ds \leq \mE\!\int_0^T\!
\!\int_G\int_{S^{d-1}}\! \theta^2 (l_t  +
U\cd\nabla l)^2 z^2 dS^{d-1}dxdt +
\!\frac{1}{2}\mE\!\int_0^T\!\!
\int_G\!\int_{S^{d-1}}\! \theta^2
(b_4^4 + 2b_4^2) z^2 dS^{d-1}dxdt \\
\ns\ds \q +  \mE\int_0^T \int_G\int_{S^{d-1}}
\theta^2 |l_t + U\cd\nabla l+ 2| Z^2
dS^{d-1}dxdt,
\end{array}
$$
by means of the inequality \eqref{bobeq3}, we
find
\begin{equation}\label{bobeq4}
\begin{array}{lll}\ds
\q \l \mE\int_G\!\int_{S^{d-1}}\!\!\!
(cT\!\!-\!2U\cd x) \theta^2(T) z^2(T)dS^{d-1}dx
\!+\!\l \!\int_G\!\int_{S^{d-1}} \!\!\!(cT \!+\! 2U\cd x)\theta^2(0) z^2(0)dS^{d-1}dx\\
\ns\ds\q  + 2(1 \!-\! c)\l\mE\!\int_0^T\!
\int_G\!\int_{S^{d-1}}\!\! \theta^2 z^2
dS^{d-1}dxdt
\!- \!3 \mE\!\int_0^T\! \int_G\!\int_{S^{d-1}}\!\! \theta^2 (b_1^2\! +\! b_4^4\!+\!b_4^2)z^2dS^{d-1}dxdt\\
\ns\ds \q - 3
|b_2|^2_{L^\infty_\cF(\O;L^\infty(0,T;C(\overline
G\times S^{d-1}\times S^{d-1})))}\mE\int_0^T
\int_G\int_{S^{d-1}} \theta^2 z^2dS^{d-1}dxdt \\
\ns\ds \leq 3\mE\int_0^T \int_G\int_{S^{d-1}}
\theta^2 \big(b_3^2  + |2  + \l x  -  c\l
t|\big)
Z^2 dS^{d-1}dxdt  \\
\ns\ds \q  - \l \mE\int_0^T
 \int_{\G_{S}^-} U\cd\nu \big[c(T-2t) - 2 U\cd
x\big]\theta^2 z^2 d\G_{S}^- dt.
\end{array}
\end{equation}
Noting that $|x|<2R$, we know that
\begin{equation}\label{bobeq5}
\left\{
\begin{array}{ll}\ds
(cT\!-\!2R)\mE\int_G\!\int_{S^{d-1}}\!\!
\theta^2(T) z^2(T)dS^{d-1}dx\leq
\mE\int_G\int_{S^{d-1}} \theta^2(T)
(cT-U\cd x) z^2(T)dS^{d-1}dx,\\
\ns\ds (cT\!-\!2R) \int_G\!\int_{S^{d-1}}\!\!
\theta^2(0) z^2(0)dS^{d-1}dx\leq
 \int_G\int_{S^{d-1}} \theta^2(0) (cT+U\cd x)
z^2(0)dS^{d-1}dx.
\end{array}
\right.
\end{equation}
Taking
$$
\begin{array}{ll}\ds
\l_1 =
\frac{3}{2(1-c)}\big(|b_1|^2_{L_{\cF}^{\infty}(0,T;L^{\infty}(G\times
S^{d-1}))} +
|b_2|^2_{L^\infty_\cF(\O;L^\infty(0,T;C(\overline
G\times S^{d-1}\times S^{d-1})))}\\
\ns\ds \qq +
|b_4|^4_{L_{\cF}^{\infty}(0,T;L^{\infty}(G\times
S^{d-1}))} +
|b_4|^2_{L_{\cF}^{\infty}(0,T;L^{\infty}(G\times
S^{d-1}))} \big),
\end{array}
$$
for any $\l\geq \l_1$, we conclude that
\begin{equation}\label{bobeq7}
\begin{array}{ll}\ds
\q 3 \mE\int_0^T \int_G \int_{S^{d-1}}\theta^2
(b_1^2 + b_4^4+b_4^2)z^2dS^{d-1}dxdt \\
\ns\ds \q +  3
|b_2|^2_{L^\infty_\cF(\O;L^\infty(0,T;C(\overline
G\times S^{d-1}\times S^{d-1})))}\mE\int_0^T
\int_G\int_{S^{d-1}} \theta^2 z^2dS^{d-1}dxdt\\
\ns\ds  \leq 2(1 - c)\l\mE\int_0^T
\int_G\int_{S^{d-1}} \theta^2 z^2 dS^{d-1}dxdt.
\end{array}
\end{equation}

From \eqref{bobeq4}--\eqref{bobeq7},
and noting that $cT>2R$, we find that
\begin{equation}\label{bobeq8}
\begin{array}{ll}\ds
\q  \mE\int_G\int_{S^{d-1}} \theta^2(T,x)
z^2(T,x)dS^{d-1}dx\\
\ns\ds \leq  C\mE\int_0^T \int_G\int_{S^{d-1}}
\theta^2 \big(b_3^2 +  2+ |\l x -  c\l t|\big)
Z^2dS^{d-1}dxdt\\
\ns\ds \q -  C\mE\int_0^T
 \int_{\G_{S}^-}  U\cd\nu
\big[c(T-2t) - 2 U\cd x\big]\theta^2 z^2
d\G_{S}^- dt.
\end{array}
\end{equation}
By the definition of $\theta$, we have
$$
e^{-c\l T^2 }\leq\theta\leq e^{4\l R^2}.
$$
This, together with the inequality
\eqref{bobeq8}, indicates that
\begin{equation}\label{bobeq9}
\begin{array}{ll}
\ds \q e^{-2c\l T^2 } \mE\int_G\int_{S^{d-1}}
z^2_T dS^{d-1}dx \\
\ns\ds \leq Ce^{8\l R^2} \left\{\mE\int_0^T
\int_G\int_{S^{d-1}}  Z^2dS^{d-1}dxdt - \mE
\int_0^T\int_{\G_{S}^-} U\cd\nu z^2 d\G_{S}^-
dt\right\},
\end{array}
\end{equation}
which implies that
\begin{equation}\label{bobeq10}
\begin{array}{lll}\ds
\q\mE\int_G\int_{S^{d-1}} z_T^2 dS^{d-1}dx \\
\ns\ds \leq Ce^{2c\l_1 T^2 + 8\l_1 R^2}
\Big\{\mE\int_0^T \int_G\int_{S^{d-1}}
Z^2dS^{d-1}dxdt - \mE \int_0^T\int_{\G_{S}^-}
U\cd\nu z^2 d\G_{S}^-
dt\Big\}\\
\ns\ds \leq e^{Cr_2^2}\Big\{\mE\int_0^T
\int_G\int_{S^{d-1}} Z^2dS^{d-1}dxdt - \mE
\int_0^T\int_{\G_{S}^-} U\cd\nu  z^2 d\G_{S}^-
dt\Big\}.
\end{array}
\end{equation}
This completes the proof.
\endpf
%%%%%%%%%%%%%%%%%%%%%%%%%%%%%%%%%%%%%%%%%%%%%%%%%%%%%%%%%%%%%%%%%%%%

\section{Proof of Theorem \ref{exact th}}

This section is addressed to a proof of Theorem
\ref{exact th}.

\vspace{0.2cm}

{\it Proof of Theorem \ref{exact th}}\,: Since
the system \eqref{csystem1} is linear, we only
need to show that the attainable set $\dbA_T$ at
time $T$ with initial datum $y(0)=0$  is
$L^2(\O,\cF_T,P;L^2(G\times S^{d-1}))$, that is,
for any $y_1\in L^2(\O,\cF_T,P;L^2(G\times
S^{d-1}))$, we can find a pair of control
$$
(u,v)\in L_{\cF}^2(0,T;L^2_{w}(\G_{S}^-))\times
L_{\cF}^2(0,T;L^2(G\times S^{d-1}))
$$
such that the solution to the system
\eqref{csystem1} satisfies that $y(T)=y_1$ in
$L^2(G\times S^{d-1})$, $P$-a.s. We achieve this
goal by the duality argument.

Let $b_1 = -a_1$, $b_2 = -a_2$, $b_3 = -a_3$ and
$b_4=0$ in the equation \eqref{csystem2}. We
introduce the following linear subspace of
$L_{\cF}^2(0,T;L^2_{w}(\G_{S}^-))\times
L_{\cF}^2(0,T;L^2(G\times S^{d-1}))$:
$$
\begin{array}{ll}\ds
\cY\=\Big\{\big(-z|_{\G_{S}^-},
Z\big)\;\Big|\;(z,Z)\hb{ solves the equation
 }\eqref{csystem2}\mbox{ with some } \\
\ns\ds \hspace{3.9cm} z_T\in
L^2(\O,\cF_T,P;L^2(G\times S^{d-1}))\Big\}
\end{array}
$$
and define a linear functional $\cL$ on $\cY$ as
follows:
$$
\cL(-z|_{\G_{S}^-},
Z)=\mathbb{E}\int_G\int_{S^{d-1}}
y_1z_TdS^{d-1}dx -
\mE\int_0^T\int_G\int_{S^{d-1}} zf dS^{d-1}dxdt.
$$
From Theorem \ref{th obser}, we see that $\cL$
is a bounded linear functional on $\cY$. By
means of the Hahn-Banach theorem, $\cL$ can be
extended to be a bounded linear functional on
the space
$L_{\cF}^2(0,T;L^2_{w}(\G_{S}^-))\times
L_{\cF}^2(0,T;L^2(G\times S^{d-1}))$. For
simplicity, we still use $\cL$ to denote this
extension. Now, by the Riesz representation
theorem,  there is a pair of random fields
$$
(u,
v)\in L_{\cF}^2(0,T;L^2_{w}(\G_{S}^-))\times
L_{\cF}^2(0,T;L^2(G\times S^{d-1}))
$$
so that
\begin{equation}\label{con eq0}
\begin{array}{ll}\ds
\q\mE\int_G\int_{S^{d-1}} y_1z_T dS^{d-1}dx -
\mE\int_0^T\int_G\int_{S^{d-1}} zfdS^{d-1}dxdt \\
\ns\ds = -\mE\int_0^T\int_{\G_{S}^-} U\cd\nu zu
d\G_{S}^- dt + \mE\int_0^T\int_G \int_{S^{d-1}}
v Z dS^{d-1}dxdt.
\end{array}
\end{equation}

We claim that this pair of random fields $(u,
v)$ is the desired controls. Indeed, by
It\^{o}'s formula, we have
\begin{equation}\label{con eq1}
\begin{array}{ll}\ds
\q\mE\int_G\int_{S^{d-1}} y(T,\cd)z_T dS^{d-1}dx\\
\ns\ds =
\mE\int_0^T\!\!\int_G\!\int_{S^{d-1}}\!\! (-z
U\cd\nabla y +a_1 yz + fz)dS^{d-1}dxdt +
\mE\int_0^T\!\!\int_G\!\int_{S^{d-1}}\!\!(a_3yZ
+ vZ
)dS^{d-1}dxdt \\
\ns\ds \q + \mE\int_0^T \int_G \int_{S^{d-1}} \(
\int_{S^{d-1}} a_2 ydS^{d-1}(V)\)z
dS^{d-1}(U)dxdt\\
\ns\ds \q
 + \mE\int_0^T \int_G \int_{S^{d-1}}
(-U\cd\nabla
z y  -  a_1yz - a_3 yZ)dS^{d-1}dxdt\\
\ns\ds \q - \mE\int_0^T \int_G \int_{S^{d-1}} \(
\int_{S^{d-1}}a_2 ydS^{d-1}(V)\)z
dS^{d-1}(U)dxdt.
\end{array}
\end{equation}
Hence,
\begin{equation}\label{con eq2}
\begin{array}{ll}\ds
\q\mE\int_G\int_{S^{d-1}} y(T,\cd)z_TdS^{d-1} dx
- \mE\int_0^T\int_G\int_{S^{d-1}}
zf dS^{d-1}dxdt \\
\ns\ds = -\mE\int_0^T\int_{\G_{S}^-}U\cd\nu zu
d\G_{S}^- dt + \mE\int_0^T\int_G\int_{S^{d-1}} v
ZdS^{d-1} dxdt.
\end{array}
\end{equation}
From   \eqref{con eq0} and  \eqref{con eq2}, we
see that
\begin{equation}\label{con eq3}
\mE\int_G\int_{S^{d-1}} y_1 z_T dS^{d-1}dx =
\mE\int_G\int_{S^{d-1}} y(T,\cd)z_T dS^{d-1}dx.
\end{equation}
Since $z_T$ can be an arbitrary element in
$L^2(\O,\cF_T,P;L^2(G\t S^{d-1}))$, from the
equality \eqref{con eq3}, we conclude that
$y(T)=y_1$ in $L^2(G\t S^{d-1})$, $P$-a.s. This
completes the proof of Theorem \ref{exact th}.
\endpf

%%%%%%%%%%%%%%%%%%%%%%%%%%%%%%%%%%%%%%%%%%%%%%%%%%%%%%%

\section{Proof of the lack of exact controllability}

The purpose of this section is to give proofs of
Theorems \ref{non control th}--\ref{non control
th2}. In order to present the key idea in the
simplest way, we only consider a very special
case of the system \eqref{csystem1}, that is,
$G=(0,1)$, $a_1=0$, $a_2=0$, $a_3 = 1$  and
$f=0$. The argument for the general case is very
similar.

\vspace{0.2cm}

{\it Proof of Theorem \ref{non control th}}\,:
The case that $v = 0$ in $L^2_\cF(0,T;L^2
(0,1))$ is covered in Theorem \eqref{non control
th1}. Hence, we only prove Theorem \ref{non
control th} for $u\equiv 0$. In this case, the
system \eqref{csystem1} reads as
\begin{equation}\label{csystem1.311}
\left\{
\begin{array}{ll}\ds
dy + y_xdt =   (y+v)dB(t)   &\mbox{ in } (0,T)\t (0,1),\\
\ns\ds y(t,0) = 0 & \mbox{ on } (0,T)\t\{0\},\\
\ns\ds y(0) = y_0 &\mbox{ in } (0,1).
\end{array}
\right.
\end{equation}
Since the system \eqref{csystem1.311} is linear,
we only need to show that the attainable set
$\dbA_T$ of this system at time $T$ for the
initial datum $y_0 =0$ is not
$L^2(\O,\cF_T,P;L^2(0,1))$. For $y_0=0$, the
solution of this system  is
\begin{equation}\label{sol2}
y(T) = \int_0^T S(T-s)\big[y(s)+v(s)\big]dB(s).
\end{equation}
Here $\{S(t)\}_{t\geq 0}$ is the semigroup
introduced in Section 2. We refer to
\cite[Chapter 6]{Prato} for the details of
establishing \eqref{sol2}. From  \eqref{sol2},
we find that $\mE (y(T))=0$. Thus, if we choose
a $y_1 \in L^2(\O,\cF_T,P;L^2(0,1))$ such that
$\mE (y_1) \neq 0$, then $y_1$ is not in
$\dbA_T$, which completes the proof.
\endpf

\vspace{0.2cm}

To prove Theorems \ref{non control
th1}--\ref{non control th2}, we first recall the
following known result.

Set
$$
\eta(t)= \left\{
\begin{array}{ll}\ds
1, &\mbox{ if } t\in
\big[(1-2^{-2i})T,(1-2^{-2i-1})T\big),\q
i=0,1,\cds,\\
\ns\ds -1,&\mbox{ otherwise  in }[0,T]
\end{array}
\right.
$$
and
\begin{equation}\label{xi}
\xi = \int_0^T \eta(t)dB(t).
\end{equation}
We have the following result.

\begin{lemma}\cite[Lemma 2.1]{Peng}\label{lm1}
It is impossible to find
$$(\varrho_1,\varrho_2)\in
L^2_\cF(0,T;\dbR)\times
C_\cF([0,T];L^2(\O;\dbR))$$ and $x\in \dbR$ such
that
\begin{equation}
\xi = x + \int_0^T \varrho_1(t)dt + \int_0^T
\varrho_2(t)dB(t).
\end{equation}
\end{lemma}

{\it Proof of Theorem \ref{non control th1}}\,:
Put
$$
\cV\=\{ v\in L^2_\cF(0,T;L^2(0,1)):\,v=0 \mbox{
in } (0,T)\times G_0\}.
$$
Let $\xi $ be given in \eqref{xi}. Choose a
$\psi\in C_0^\infty(G_0)$ such that
$|\psi|_{L^2(G)}=1$ and set $y_T = \xi\psi$. We
will show that $y_T$ cannot be attained for any
$y_0\in \dbR$, $u\in L^2_\cF(0,T;\dbR)$ and
$v\in \cV$. This goal is achieved by the
contradiction argument. If there exist a
$y_0\in\dbR$, a $u\in L^2_\cF(0,T;\dbR)$ and a
$v\in \cV$ such that the corresponding solution
$y(\cd)$ satisfies $y(T)=y_T$, then by the
definition of the solution to \eqref{csystem1},
we obtain that
\begin{equation}\label{thn1-eq2}
\xi=\int_G  y_T\psi dx  = \int_G  y_0\psi dx +
\int_0^T \(\int_G \psi_x y dx \)dt + \int_0^T
\(\int_G \psi y dx\)dB(t).
\end{equation}
It is clear that both $ \ds\int_G \psi_x y dx $
and $\ds\int_G  \psi y dx$ belong to
$C_\cF([0,T];L^2(\O;\dbR))$. This, together with
\eqref{thn1-eq2}, contradicts Lemma \ref{lm1}.
\endpf
\medskip

{\it Proof of Theorem \ref{non control th2}}\,:
The proof is similar to the one for Theorem
\ref{non control th1}.

\vspace{0.1cm}

Let $\xi $ be given by \eqref{xi}. Choose a
$\psi\in C_0^\infty(G)$ such that
$|\psi|_{L^2(G)}=1$ and set $y_T = \xi\psi$. We
will show that $y_T$ cannot be attained for any
$y_0\in \dbR$, $u\in L^2_\cF(0,T;\dbR)$ and
$\ell\in L^2_\cF(0,T;L^2(0,1))$. It is done by
the contradiction argument too. If there exist a
$u\in L^2_\cF(0,T;\dbR)$ and an $\ell\in
L^2_\cF(0,T;L^2(0,1))$ such that the
corresponding solution $y(\cd)$ satisfies
$y(T)=y_T$, then, from the definition of the
solution to \eqref{csystem1}, we obtain
\begin{equation}\label{thn2-eq2}
\begin{array}{ll}\ds
\xi\3n &\ds=\int_G  y_T\psi dx \\
\ns&\ds = \int_G  y_0\psi dx + \int_0^T
\(\int_G \psi_x y dx + \int_G \psi \ell dx \)dt
+ \int_0^T \(\int_G \psi y dx\)dB(t).
\end{array}
\end{equation}
It is clear that $\ds\int_G \psi_x y dx + \int_G
\psi \ell dx \in L^2_\cF(0,T;\dbR)$ and
$\ds\int_G \psi y\in C_\cF([0,T];L^2(\O;\dbR))$.
These, together with \eqref{thn2-eq2},
contradict Lemma \ref{lm1}.
\endpf

%%%%%%%%%%%%%%%%%%%%%%%%%%%%%%%%%%%%%%%%%%%%%%%%%%%%%%%%%%%%%%%%%%%%%%%%%%%%%%%%%%%%

% {\bf \large
%Acknowledgments.}

\end{document}